\documentclass[a4paper,12pt]{amsart}
\usepackage{amssymb}
\usepackage{ifthen}
\usepackage{graphicx}
\nonstopmode

 \def\registered{
 {\ooalign{\hfil\raise .00ex\hbox{\scriptsize R}\hfil\crcr\mathhexbox20D}}}
\newtheorem{thm}{Theorem}
\newtheorem{cor}{Corollary}
\newtheorem{lem}{Lemma}

\newtheorem{claim}{Claim}

\newtheorem{conj}{Conjecture}
\newtheorem{prob}{Problem}

\theoremstyle{definition}

\newenvironment{rem}{%
\bigskip
\noindent \textsl{{\sl Remark. }}}{\bigskip}
\newenvironment{rems}{%
\bigskip
\noindent \textsl{{\sl Remarks. }}}{\bigskip}

\newenvironment{pf}[1][]{%
 \vskip 1mm
 \noindent
 \ifthenelse{\equal{#1}{}}%
  {{\slshape Proof. }}%
  {{\slshape #1.} }%
 }%
{\qed\medskip}

\newcounter{alphabet}
\newcounter{tmp}
\newenvironment{Thm}[1][]{\refstepcounter{alphabet}%
\bigskip%
\noindent%
{\bf Theorem \Alph{alphabet}}%
\ifthenelse{\equal{#1}{}}{}{ (#1)}%
{\bf .} \itshape}{\vskip 8pt}

\makeatletter
\newcommand{\Ref}[1]{\@ifundefined{r@#1}{}{\setcounter{tmp}{\ref{#1}}\Alph{tmp}}}
\makeatother

\newenvironment{Lem}[1][]{\refstepcounter{alphabet}%
\bigskip%
\noindent%
{\bf Lemma \Alph{alphabet}}%
{\bf .} \itshape}{\vskip 8pt}

\newcommand{\IR}{{\mathbb R}}

\newcommand{\IC}{{\mathbb C}}
\newcommand{\ID}{{\mathbb D}}

\def\be{\begin{equation}}
\def\ee{\end{equation}}

\newcommand{\bee}{\begin{enumerate}}
\newcommand{\eee}{\end{enumerate}}

\newcommand{\blem}{\begin{lem}}
\newcommand{\elem}{\end{lem}}
\newcommand{\bthm}{\begin{thm}}
\newcommand{\ethm}{\end{thm}}
\newcommand{\bcor}{\begin{cor}}
\newcommand{\ecor}{\end{cor}}
\newcommand{\beg}{\begin{examp}}
\newcommand{\eeg}{\end{examp}}
\newcommand{\begs}{\begin{examples}}
\newcommand{\eegs}{\end{examples}}
\newcommand{\bdefe}{\begin{defin}}
\newcommand{\edefe}{\end{defin}}
\newcommand{\bprob}{\begin{prob}}
\newcommand{\eprob}{\end{prob}}
\newcommand{\bques}{\begin{ques}}
\newcommand{\eques}{\end{ques}}
\newcommand{\bei}{\begin{itemize}}
\newcommand{\eei}{\end{itemize}}

\newcommand{\bde}{\begin{deter}}
\newcommand{\ede}{\end{deter}}
\newcommand{\bca}{\begin{case}}
\newcommand{\eca}{\end{case}}
\newcommand{\bcl}{\begin{claim}}
\newcommand{\ecl}{\end{claim}}
\newcommand{\bcon}{\begin{conj}}
\newcommand{\econ}{\end{conj}}
\newcommand{\bcons}{\begin{conjs}}
\newcommand{\econs}{\end{conjs}}
\newcommand{\bprop}{\begin{propo}}
\newcommand{\eprop}{\end{propo}}
\newcommand{\br}{\begin{rem}}
\newcommand{\er}{\end{rem}}
\newcommand{\brs}{\begin{rems}}
\newcommand{\ers}{\end{rems}}
\newcommand{\bo}{\begin{obser}}
\newcommand{\eo}{\end{obser}}
\newcommand{\bos}{\begin{obsers}}
\newcommand{\eos}{\end{obsers}}
\newcommand{\bpf}{\begin{pf}}
\newcommand{\epf}{\end{pf}}
\newcommand{\ba}{\begin{array}}
\newcommand{\ea}{\end{array}}
\newcommand{\beq}{\begin{eqnarray}}
\newcommand{\beqq}{\begin{eqnarray*}}
\newcommand{\eeq}{\end{eqnarray}}
\newcommand{\eeqq}{\end{eqnarray*}}

\newcommand{\ds}{\displaystyle}

\newcounter{minutes}
\divide\time by 60
\newcounter{hours}
\multiply\time by 60 \addtocounter{minutes}{-\time}

\begin{document}
\title[The minimal surfaces over the slanted
half-planes, vertical strips and single slit] {The minimal surfaces
over the slanted half-planes, vertical strips and single slit}
\thanks{
File:~\jobname .tex,
          printed: \number\day-\number\month-\number\year,
          \thehours.\ifnum\theminutes<10{0}\fi\theminutes}

\author{Liulan Li, S. Ponnusamy $^\dagger $ and M. Vuorinen
}
\address{Liulan Li, Department of Mathematics and Computational Science,
Hengyang Normal University, Hengyang,  Hunan 421008, People's
Republic of China}
\email{lanlimail2008@yahoo.com.cn}
\address{S. Ponnusamy, Department of Mathematics, Indian Institute of Technology Madras, Chennai 600036, India}
\email{samy@iitm.ac.in}
\address{M. Vuorinen, Department of Mathematics,
University of Turku,  Turku 20014, Finland.}
\email{vuorinen@utu.fi}


\keywords{Univalent Harmonic mapping, slanted half-plane mapping, slit mapping, strip mapping,
convex in the real direction, minimal surface. \\
$
^\dagger$ {\tt Corresponding author}
}



\begin{abstract}
In this paper, we discuss the minimal surfaces over
the slanted half-planes, vertical strips, and single slit whose slit
lies on the negative real axis. The representation of these minimal
surfaces and the corresponding harmonic mappings are obtained explicitly.
Finally, we illustrate the harmonic mappings of each of these cases
together with their minimal surfaces pictorially with the help of
mathematica.
\end{abstract}

\thanks{The research of the first author was supported by NSF of Hunan (No. 10JJ4005), Hunan Provincial Education
Department (No. 11B019) and partly supported by the construct
program of the key discipline in Hunan province.
}


\maketitle
\pagestyle{myheadings}
\markboth{Liulan Li, S. Ponnusamy and M. Vuorinen}{The minimal surfaces over the slanted half-planes, vertical strips and single slit}

\section{Introduction}

A planar harmonic mapping in the unit disk $\ID =\{z:\, |z|<1\}$ is a complex-valued harmonic function $f(z)$, defined
on $\ID$. The mapping $f$ has a canonical decomposition $f=h+\overline{g}$, where $h$ and $g$ are analytic on $\ID$ and
$g(0)=0$. The mapping $f$ is locally univalent in $\ID$ if and only if its Jacobian
$J_f(z)=|h'(z)|^2 -|g'(z)|^2$ does not vanish in $\ID$. It is said to be
sense-preserving on $\ID$ if and only if $J_f(z)>0$, or equivalently
if $h'(z)\neq 0$ in $\ID$ and $f$ satisfies the elliptic partial differential equation
$$\overline{f_{\overline  z}(z)}=\omega (z)f_z(z)
$$
in $\ID$, where the dilatation $\omega (z)=g'(z)/h'(z)$ has the property that $|\omega (z)|<1$ in $\ID$.

Planar univalent harmonic mappings are used in the study of the Gaussian curvature of
nonparametric minimal surfaces over simply connected domains  (see for example \cite{Du,H52}).
After the publication of landmark paper of Clunie and Sheil-Small \cite{Clunie-Small-84}, considerable
interest in the function theoretic properties of harmonic functions, quite apart from this connection, was generated.
Since then the study of  univalent harmonic mappings has gained much attention. The case where
$\omega (z)$ is a finite Blaschke product
is of special interest since this case arises in many different contexts (see \cite{HS88,SS89}). In the present
paper we shall explicitly  study the connection between certain classes of harmonic univalent mappings
and the theory of minimal surfaces.


Let $S$ be a nonparametric minimal surface over a simply connected domain $\Omega$ in $\IC$ given by
$$S=\{(u, v, F(u,v) ):\, u+iv\in\Omega\},
$$
where we have identified $\IR^2$ with the complex plane in describing the domain of $F$. The following result due to
Weierstrass-Enneper representation provides the close link between harmonic univalent mappings and
the associated minimal surfaces.
Then $S$ is a minimal surface if and only if $S$ has the
representation of the form
$$S=\left\{\left({\rm Re} \int^{z}_{0}\phi_1(t)\,dt+c_1,
{\rm Re} \int^{z}_{0}\phi_2(t)\,dt+c_2, {\rm Re} \int^{z}_{0}\phi_3(t)\,dt+c_3 \right):\, z\in\ID\right\},
$$
where $\phi_1,\ \phi_2,\ \phi_3$ are analytic in $\ID$,
\be\label{eq2}
\phi^2_1+\phi^2_2+\phi^2_3=0,\mbox{ and } f=u+iv={\rm Re}
\int^{z}_{0}\phi_1(t)\,dt+i{\rm Re} \int^{z}_{0}\phi_2(t)\,dt+c
\ee
is a sense-preserving univalent harmonic mapping from $\ID$ onto $\Omega$.
For this case, we call $S$  a minimal surface over $\Omega$ with
the projection $f=u+iv$.

Further basic information about harmonic mappings and their relation to minimal surfaces may be found in the book of
Duren \cite{Du}. For instance, the following formulation is well-known (see for instance \cite[Section 10.2]{Du}).

\begin{Thm} 
\label{Thm A}
If $f=h+\overline{g}$ is a harmonic mapping of the form \eqref{eq2} with the dilatation
$\omega =b^2$, where $b(z)=\pm z$, then we have
$$\phi_1=h'+g',\ \phi_2=-i(h'-g'),\ \phi_3=2ibh'.$$
\end{Thm}

Using this, Jun \cite{Jun06} has considered the minimal surfaces associated with the
harmonic mappings especially when $\Omega=\{w:\, {\rm Im\,} w>0\}$. His main result, which is easy to prove, will
now be recalled for the sake of convenient reference.

\begin{Thm}{\rm (\cite{Jun06})}\label{Thm B}
Let $\Omega=\{w:\, {\rm Im\,} w>0\}$ and $p=p_1+ip_2$ be a fixed
point in $\Omega$, where $p_1,\ p_2\in\IR$. If $S$ is a minimal
surface over $\Omega$ with the projection $f=h+\overline{g}$, where
$\omega (z)=\frac{g'(z)}{h'(z)}=b^2(z)=z^2$, $b(z)=\pm z$ and
$f(0)=p$, then $S=\{(u, v, F(u,v) ):\, u+iv\in\Omega\},$ where
\beqq
u&=&{\rm Re\,}f(z)=p_1+\frac{ip_2}{2}\left[\left(\frac{1}{2}\log\frac{1+z}{1-z}+\frac{z}{(1-z)^2}\right)-
\overline{\left(\frac{1}{2}\log\frac{1+z}{1-z}+\frac{z}{(1-z)^2}\right)}\right],\\
v&=&{\rm Im\,}f(z)=\frac{p_2}{2}\left[\frac{1+z}{1-z}+\overline{\left(\frac{1+z}{1-z}\right)}\right],\\
F&=&\pm p_2{\rm Re}\left(\frac{z}{(1-z)^2}- \frac{1}{2}\log\frac{1+z}{1-z}\right).
\eeqq
\end{Thm}

The class ${\mathcal S}_H$ of sense-preserving harmonic univalent mappings $f=h+\overline{g}$
(normalized so that $f(0)=0=h(0)$ and $f_z(0)=1$)
together with its many geometric subclasses have been extensively studied
(see \cite{Clunie-Small-84,Du}). Let ${\mathcal S}_H^0$ be the subset
of all $f\in {\mathcal S}_H$ in which $b_1=f_{\overline{z}}(0)=0$.
We remark that the familiar class ${\mathcal S}$ of normalized
analytic univalent functions is contained in ${\mathcal S}_H^0$.
Every  $f \in {\mathcal S}_H$ admits the complex dilatation $\omega
$ of $f$ which satisfies $|\omega (z)|<1$ in $\ID$. When $f \in {\mathcal S}_H^0$,
we also have $\omega '(0)=0$.
\medskip

In this paper, we discuss the minimal surfaces over the slanted
half-planes, vertical strips, and single slit whose slit lies on the
negative real axis. Slanted half-plane mappings are well suited in the study
of convolution of harmonic mappings (see \cite{DN}). Since the slanted half-planes and vertical
strips are convex domains, the following result of Clunie and Sheil-Small
is applicable for these cases.

\begin{Lem}{\rm \cite{Clunie-Small-84}}\label{LemA}
If $f=h+\overline{g}$ is a sense-preserving univalent
mapping such that $f(\ID)$ is a convex domain, then the function
$h+e^{i\beta}g$ is univalent for each $\beta$, $0\leq\beta<2\pi$.
\end{Lem}

\section{Slanted half-plane mappings}\label{half-plane}

Throughout this section, we let $H_\gamma :=\{w:\,{\rm
Re\,}(e^{i\gamma}w) >-1/2\}$ be a slanted half-plane with the parameter $\gamma$,
where $0\leq\gamma<2\pi$.

\begin{thm}\label{1}
Let $S$ be a minimal surface over $H_\gamma$ with the projection
$f=h+\overline{g}$, whose dilatation $\omega=g'/h'=b^2$, where
$b(z)=\pm z$. Then
$$S=\{(u, v, F(u,v) ):\, u+iv\in H_\gamma\} = \{(u(z), v(z), F(u(z),v(z)) ):\, z\in \ID\},
$$
where
\beqq
u& =& \frac{\pi \sin\gamma}{4}-\cos\gamma-{\rm Im}\left(\frac{\sin\gamma}{4}\log\frac{z-e^{-i\gamma}}{z+e^{-i\gamma}}
 +\frac{\sin2\gamma}{4(z-e^{-i\gamma})}\right) \\
&& \hspace{.3cm}  -{\rm Re}\left(\frac{\cos\gamma}{2(z-e^{-i\gamma})^2}+\frac{3}{4(z-e^{-i\gamma})}\right),\\
v&=& \frac{\pi \cos\gamma}{4}+\sin\gamma-{\rm Im}\left(\frac{\cos\gamma}{4}\log\frac{z-e^{-i\gamma}}{z+e^{-i\gamma}}
-\frac{3}{4(z-e^{-i\gamma})}\right)\\
&& \hspace{.3cm} +{\rm Re}\left(\frac{\sin2\gamma}{4(z-e^{-i\gamma})}-\frac{\sin\gamma}{2(z-e^{-i\gamma})^2}\right),\\
F&=&\pm{\rm Re}\left[\frac{\sin2\gamma}{4}\log\frac{z+e^{-i\gamma}}{z-e^{-i\gamma}}
+\frac{1}{2}e^{i(\gamma+\frac{\pi}{2})}\frac{1}{z-e^{-i\gamma}}
+\frac{i}{2}\frac{1}{(z-e^{-i\gamma})^2}\right]+c,
\eeqq
if $\ds \gamma\in\Big \{\frac{\pi}{4}, \frac{3\pi}{4}, \frac{5\pi}{4}, \frac{7\pi}{4}\Big \}$;
\beqq
u& =& {\rm Im}\left(\frac{\sin\gamma}{2(1-\sin
2\gamma)}\log\frac{z+ie^{i\gamma}}{ie^{i\gamma}}+\frac{\sin\gamma}{2(1+\sin
2\gamma)}\log\frac{z-ie^{i\gamma}}{-ie^{i\gamma}}-\frac{\sin\gamma}{\cos^2
2\gamma}\log\frac{z-e^{-i\gamma}}{-e^{-i\gamma}}\right)\\
&& \hspace{.3cm}  -\frac{\cos\gamma}{\cos2\gamma}-\frac{\cos\gamma}{\cos2\gamma}{\rm
Re\,}\frac{e^{-i\gamma}}{z-e^{-i\gamma}},\\
v&=&{\rm Im}\left(\frac{\cos\gamma}{2(1-\sin
2\gamma)}\log\frac{z+ie^{i\gamma}}{ie^{i\gamma}}+\frac{\cos\gamma}{2(1+\sin
2\gamma)}\log\frac{z-ie^{i\gamma}}{-ie^{i\gamma}}-\frac{\cos\gamma}{\cos^2
2\gamma}\log\frac{z-e^{-i\gamma}}{-e^{-i\gamma}}\right)\\
&& \hspace{.3cm} -\frac{\sin\gamma}{\cos2\gamma}-\frac{\sin\gamma}{\cos2\gamma}{\rm
Re\,}\frac{e^{-i\gamma}}{z-e^{-i\gamma}},\\
F&=& \pm{\rm  Re}\left[\frac{\log(z+ie^{i\gamma})}{2(1-\sin
2\gamma)}-\frac{\log(z-ie^{i\gamma})}{2(1+\sin 2\gamma)}-\frac{\sin
2\gamma}{\cos^2
2\gamma}\log(z-e^{-i\gamma})-\frac{ie^{-i\gamma}}{(z-e^{-i\gamma})\cos
2\gamma}\right]+c,
\eeqq
if $\ds \gamma\notin\Big\{\frac{\pi}{4}, \frac{3\pi}{4}, \frac{5\pi}{4}, \frac{7\pi}{4}\Big\}$.
\end{thm}\bpf
Let $f=h+\overline{g}\in {\mathcal S}_H^0$ and $f(\ID)=H_{\gamma}$. Then, we have
$${\rm Re\,}(e^{i\gamma}f(z)) ={\rm Re\,}[e^{i\gamma}(h(z)+e^{-2i\gamma}g(z))]>-\frac{1}{2}, \quad z\in\ID,
$$
so that $(h+e^{-2i\gamma}g)(\ID)=H_{\gamma}$ and by Lemma \Ref{LemA},  $h+e^{-2i\gamma}g$ is conformal (univalent) mapping
from $\ID$ onto $H_\gamma$.

We now consider the function $h+e^{-2i\gamma}g$. We may conveniently normalize it in such a way
that $f(0)=h(0)=g(0)=0$. Then $h(0)+e^{-2i\gamma}g(0)=0$. We further
assume that
$$h(e^{-i\gamma})+e^{-2i\gamma}g(e^{-i\gamma})=\infty~\mbox{ and }~ h(e^{-i(\pi +\gamma)})+e^{-2i\gamma}g(e^{-i(\pi +\gamma)})=-\frac{1}{2}e^{-i\gamma}.
$$
By the uniqueness of the Riemann mapping theorem, these observations led to the representation (see also \cite[Lemma 1]{DN})
\be\label{eq4}
h(z)+e^{-2i\gamma}g(z)=\frac{z}{1-e^{i\gamma}z}
\ee
from which we obtain
\be\label{eq4a}
g(z)=-\frac{1}{z-e^{-i\gamma}} -e^{2i\gamma}h(z)- e^{i\gamma}
\ee
and
$$h'(z)+e^{-2i\gamma}g'(z)=\frac{1}{(1-e^{i\gamma}z)^2}.
$$
Solving this together with $g'(z)=z^2h'(z)$ gives
$$h'(z)=\frac{1}{(z^2+e^{2i\gamma})(z-e^{-i\gamma})^2}~\mbox{ and }~
g'(z)=\frac{z^2}{(z^2+e^{2i\gamma})(z-e^{-i\gamma})^2}.
$$
It is convenient to write $h'(z)$ in the form
\be\label{eq4b}
h'(z)= \frac{1}{(z-e^{i(\gamma +\pi/2)})(z-e^{i(\gamma -\pi/2)})(z-e^{-i\gamma})^2}.
\ee
In order to determine $h(z)$ explicitly, we need to decompose it into partial fractions, and it is also clear that
we need to deal with the cases where
$$ \gamma\in\left \{\frac{\pi}{4}, \frac{3\pi}{4}, \frac{5\pi}{4}, \frac{7\pi}{4}\right \} \mbox{ and }~
\gamma\notin\left \{\frac{\pi}{4}, \frac{3\pi}{4}, \frac{5\pi}{4}, \frac{7\pi}{4}\right \}.
$$
 \medskip

\noindent {\bf Case 1:}\quad Let $\gamma=\frac{\pi}{4}.$

In this case, $h'(z)$ given by \eqref{eq4b} takes the form
$$h'(z)=\frac{1}{(z+e^{-\frac{i\pi}{4}})(z-e^{-\frac{i\pi}{4}})^3}
$$
so that $h'(z)$ has a simple pole at $z=-e^{-\frac{i\pi}{4}}$ and a pole of order $3$ at $z=e^{-\frac{i\pi}{4}}$. We see that
$$h'(z)=\frac{i}{8}e^{\frac{i\pi}{4}}\left (\frac{1}{z-e^{-\frac{i\pi}{4}}}-\frac{1}{z+e^{-\frac{i\pi}{4}}}\right )
-\frac{i}{4}\frac{1}{(z-e^{-\frac{i\pi}{4}})^2}+\frac{1}{2}e^{\frac{i\pi}{4}}\frac{1}{(z-e^{-\frac{i\pi}{4}})^3},
$$
Integration from $0$ to $z$ gives
\be\label{eq4f}
h(z)=\left[\frac{1}{8}e^{\frac{3i\pi}{4}}\log\frac{z-e^{-\frac{i\pi}{4}}}{z+e^{-\frac{i\pi}{4}}}+
\frac{i}{4}\frac{1}{z-e^{-\frac{i\pi}{4}}}-\frac{1}{4}e^{\frac{i\pi}{4}}\frac{1}{(z-e^{-\frac{i\pi}{4}})^2}\right]
-\frac{1}{2}e^{-\frac{i\pi}{4}}+\frac{\pi}{8}e^{\frac{i\pi}{4}}.
\ee
Equation \eqref{eq4a} for $\gamma=\frac{\pi}{4}$ gives
$$g(z)=-\frac{1}{z-e^{-\frac{i\pi}{4}}} -ih(z)- e^{\frac{i\pi}{4}}
$$
so that
$$h(z)+g(z)=-\frac{1}{z-e^{-\frac{i\pi}{4}}} +\sqrt{2}e^{-\frac{i\pi}{4}}h(z)- e^{\frac{i\pi}{4}}
$$
and thus, substituting the expression for $h(z)$ defined by \eqref{eq4f} yields that
$$h(z)+g(z)=\frac{\sqrt{2}\pi}{8}-\frac{\sqrt{2}}{2}+\frac{i\sqrt{2}}{8}\log\frac{z-e^{-\frac{i\pi}{4}}}{z+e^{-\frac{i\pi}{4}}} -\frac{3-i}{4}\frac{1}{z-e^{-\frac{i\pi}{4}}}
-\frac{\sqrt{2}}{4}\frac{1}{(z-e^{-\frac{i\pi}{4}})^2}
$$
and similarly
$$h(z)-g(z)=\frac{i\sqrt{2}\pi}{8}+\frac{i\sqrt{2}}{2}-\frac{\sqrt{2}}{8}\log\frac{z-e^{-\frac{i\pi}{4}}}{z+e^{-\frac{i\pi}{4}}}
+\frac{3+i}{4}\frac{1}{z-e^{-\frac{i\pi}{4}}}
-\frac{i\sqrt{2}}{4}\frac{1}{(z-e^{-\frac{i\pi}{4}})^2}.
$$
As $u={\rm Re\,}f(z)={\rm Re\,}(h(z)+g(z))$ and $v={\rm Im\,}f(z)={\rm Im\,}(h(z)-g(z))$, the last two equalities give
$$u=\frac{\sqrt{2}\pi}{8}-\frac{\sqrt{2}}{2}-{\rm
Im}\left(\frac{\sqrt{2}}{8}\log\frac{z-e^{-\frac{i\pi}{4}}}{z+e^{-\frac{i\pi}{4}}}
+\frac{1}{4}\frac{1}{z-e^{-\frac{i\pi}{4}}}\right)-{\rm
Re}\left(\frac{\sqrt{2}}{4}\frac{1}{(z-e^{-\frac{i\pi}{4}})^2}+\frac{3}{4}\frac{1}{z-e^{-\frac{i\pi}{4}}}\right),
$$
and
$$v=\frac{\sqrt{2}\pi}{8}+\frac{\sqrt{2}}{2}-{\rm
Im}\left(\frac{\sqrt{2}}{8}\log\frac{z-e^{-\frac{i\pi}{4}}}{z+e^{-\frac{i\pi}{4}}}-\frac{3}{4}\frac{1}{z-e^{-\frac{i\pi}{4}}}\right)+{\rm
Re}\left(\frac{1}{4}\frac{1}{z-e^{-\frac{i\pi}{4}}}-\frac{\sqrt{2}}{4}\frac{1}{(z-e^{-\frac{i\pi}{4}})^2}\right).
$$

Finally, as $b(z)=\pm z$, Theorem \Ref{Thm A} gives,
\beqq
\phi_3(z) &=& 2ibh'(z)=\pm 2i\frac{z}{(z+e^{-\frac{i\pi}{4}})(z-e^{-\frac{i\pi}{4}})^3}\\
&=&\pm 2i\left[\frac{i}{8}\frac{1}{z+e^{-\frac{i\pi}{4}}}-\frac{i}{8}\frac{1}{z-e^{-\frac{i\pi}{4}}}
+\frac{1}{4}e^{\frac{i\pi}{4}}\frac{1}{(z-e^{-\frac{i\pi}{4}})^2}+\frac{1}{2}\frac{1}{(z-e^{-\frac{i\pi}{4}})^3}\right].
\eeqq
and therefore,
\beqq
F(z)&=&{\rm Re\,} \int^{z}_{0}\phi_3(z)\,dz+c\\
&=&\mp {\rm
Re}\left[\frac{1}{4}\log\frac{z+e^{-\frac{i\pi}{4}}}{z-e^{-\frac{i\pi}{4}}}+\frac{1}{2}e^{\frac{3i\pi}{4}}\frac{1}{z-e^{-\frac{i\pi}{4}}}
+\frac{i}{2}\frac{1}{(z-e^{-\frac{i\pi}{4}})^2}\right]+c.
\eeqq

\medskip

\noindent {\bf Case 2:}\quad Let $\gamma=\frac{3\pi}{4}$.

In this case, $h'(z)$ given by \eqref{eq4b} takes the form
$$h'(z)=\frac{1}{(z-e^{\frac{i\pi}{4}})(z+e^{\frac{i\pi}{4}})^3}
$$
and the partial fraction expansion gives
$$h'(z)=\frac{1}{8}e^{\frac{i\pi}{4}}\left (\frac{1}{z+e^{\frac{i\pi}{4}}}-\frac{1}{z-e^{\frac{i\pi}{4}}}\right )
+\frac{i}{4}\frac{1}{(z+e^{\frac{i\pi}{4}})^2}-\frac{1}{2}e^{-\frac{i\pi}{4}}\frac{1}{(z+e^{\frac{i\pi}{4}})^3}.
$$
Integration from $0$ to $z$ gives
\be\label{eq4g}
 h(z)=\frac{1}{8}e^{\frac{i\pi}{4}}\log\frac{z+e^{\frac{i\pi}{4}}}{z-e^{\frac{i\pi}{4}}}-\frac{i}{4}\frac{1}{z+e^{\frac{i\pi}{4}}}
+\frac{1}{4}e^{-\frac{i\pi}{4}}\frac{1}{(z+e^{\frac{i\pi}{4}})^2}+\frac{1}{2}e^{\frac{i\pi}{4}}-\frac{i\pi}{8}e^{\frac{i\pi}{4}}.
\ee
Using \eqref{eq4a} for $\gamma=\frac{3\pi}{4}$, we see that
$$h(z)+g(z)=-\frac{1}{z+e^{\frac{i\pi}{4}}} +\sqrt{2}e^{\frac{i\pi}{4}}h(z)+ e^{-\frac{i\pi}{4}}
$$
and
$$h(z)-g(z)=\frac{1}{z+e^{\frac{i\pi}{4}}} +\sqrt{2}e^{-\frac{i\pi}{4}}h(z)- e^{-\frac{i\pi}{4}}
$$
which, by \eqref{eq4g}, simplify to
$$h(z)+g(z)=\frac{\sqrt{2}\pi}{8}+\frac{\sqrt{2}}{2}+\frac{i\sqrt{2}}{8}\log\frac{z+e^{\frac{i\pi}{4}}}{z-e^{\frac{i\pi}{4}}}
-\frac{3+i}{4}\frac{1}{z+e^{\frac{i\pi}{4}}}+\frac{\sqrt{2}}{4}\frac{1}{(z+e^{\frac{i\pi}{4}})^2}
$$
and
$$h(z)-g(z)=-\frac{i\sqrt{2}\pi}{8}+\frac{i\sqrt{2}}{2}+\frac{\sqrt{2}}{8}\log\frac{z+e^{\frac{i\pi}{4}}}{z-e^{\frac{i\pi}{4}}}
+\frac{3-i}{4}\frac{1}{z+e^{\frac{i\pi}{4}}}-\frac{i\sqrt{2}}{4}\frac{1}{(z+e^{\frac{i\pi}{4}})^2},
$$
respectively. As $u={\rm Re\,}f(z)={\rm Re\,}(h(z)+g(z))$ and $v={\rm Im\,}f(z)={\rm Im\,}(h(z)-g(z))$, it follows easily that
$$u=\frac{\sqrt{2}\pi}{8}+\frac{\sqrt{2}}{2}-{\rm
Im}\left(\frac{\sqrt{2}}{8}\log\frac{z+e^{\frac{i\pi}{4}}}{z-e^{\frac{i\pi}{4}}}
-\frac{1}{4}\frac{1}{z+e^{\frac{i\pi}{4}}}\right)+{\rm
Re}\left(\frac{\sqrt{2}}{4}\frac{1}{(z+e^{\frac{i\pi}{4}})^2}-\frac{3}{4}\frac{1}{z+e^{\frac{i\pi}{4}}}\right),
$$
and
$$v=-\frac{\sqrt{2}\pi}{8}+\frac{\sqrt{2}}{2}+{\rm
Im}\left(\frac{\sqrt{2}}{8}\log\frac{z+e^{\frac{i\pi}{4}}}{z-e^{\frac{i\pi}{4}}}
+\frac{3}{4}\frac{1}{z+e^{\frac{i\pi}{4}}}\right)-{\rm
Re}\left(\frac{\sqrt{2}}{4}\frac{1}{(z+e^{\frac{i\pi}{4}})^2}+\frac{1}{4}\frac{1}{z+e^{\frac{i\pi}{4}}}\right).
$$
Moreover, in this case Theorem \Ref{Thm A} implies that
\beqq
\phi_3(z) &=&  2ibh'(z)= \pm 2i\frac{z}{(z-e^{\frac{i\pi}{4}})(z+e^{\frac{i\pi}{4}})^3}\\
&= &\pm 2i\left[-\frac{i}{8}\frac{1}{z-e^{\frac{i\pi}{4}}}+\frac{i}{8}\frac{1}{z+e^{\frac{i\pi}{4}}}
+\frac{1}{4}e^{\frac{3i\pi}{4}}\frac{1}{(z+e^{\frac{i\pi}{4}})^2}+\frac{1}{2}\frac{1}{(z+e^{\frac{i\pi}{4}})^3}\right].
\eeqq
Integration from $0$ to $z$ gives
$$F(z)=\mp{\rm
Re}\left[-\frac{1}{4}\log\frac{z-e^{\frac{i\pi}{4}}}{z+e^{\frac{i\pi}{4}}}-\frac{1}{2}e^{\frac{i\pi}{4}}\frac{1}{z+e^{\frac{i\pi}{4}}}
+\frac{i}{2}\frac{1}{(z+e^{\frac{i\pi}{4}})^2}\right]+c.
$$

\medskip

\noindent {\bf Case 3:}\quad Let  $\gamma=\frac{5\pi}{4}$.

In this case, $h'(z)$ given by \eqref{eq4b} takes the form
$$h'(z)=\frac{1}{(z-e^{-\frac{i\pi}{4}})(z+e^{-\frac{i\pi}{4}})^3}
$$
and therefore,
$$ h'(z)=\frac{1}{8}e^{-\frac{i\pi}{4}} \left (\frac{1}{z+e^{-\frac{i\pi}{4}}}-\frac{1}{z-e^{-\frac{i\pi}{4}}}\right )
-\frac{1}{4}e^{\frac{i\pi}{2}}\frac{1}{(z+e^{-\frac{i\pi}{4}})^2}-\frac{1}{2}e^{\frac{i\pi}{4}}\frac{1}{(z+e^{-\frac{i\pi}{4}})^3}.
$$
Integration from $0$ to $z$ gives
\be\label{eq4h}
 h(z)=\frac{1}{8}e^{-\frac{i\pi}{4}}\log\frac{z+e^{-\frac{i\pi}{4}}}{z-e^{-\frac{i\pi}{4}}}+\frac{i}{4}\frac{1}{z+e^{-\frac{i\pi}{4}}}
+\frac{1}{4}e^{\frac{i\pi}{4}}\frac{1}{(z+e^{-\frac{i\pi}{4}})^2}-\frac{\pi}{8}e^{\frac{i\pi}{4}}+\frac{1}{2}e^{-\frac{i\pi}{4}}.
\ee
Using \eqref{eq4a} for $\gamma=\frac{5\pi}{4}$, we see that
$$h(z)+g(z)=-\frac{1}{z+e^{-\frac{i\pi}{4}}} +\sqrt{2}e^{-\frac{i\pi}{4}}h(z)+ e^{\frac{i\pi}{4}}
$$
and
$$h(z)-g(z)=\frac{1}{z+e^{-\frac{i\pi}{4}}} +\sqrt{2}e^{\frac{i\pi}{4}}h(z)- e^{\frac{i\pi}{4}}.
$$
As in the earlier two cases, a routine computation with the help of \eqref{eq4h} shows that
$$u={\rm
Im}\left(\frac{\sqrt{2}}{8}\log\frac{z+e^{-\frac{i\pi}{4}}}{z-e^{-\frac{i\pi}{4}}}
-\frac{1}{4}\frac{1}{z+e^{-\frac{i\pi}{4}}}\right)
+{\rm Re}\left(\frac{\sqrt{2}}{4}\frac{1}{(z+e^{-\frac{i\pi}{4}})^2}-\frac{3}{4}\frac{1}{z+e^{-\frac{i\pi}{4}}}\right)
-\frac{\sqrt{2}\pi}{8}+\frac{\sqrt{2}}{2},
$$
and
$$v={\rm
Im}\left(\frac{\sqrt{2}}{8}\log\frac{z+e^{-\frac{i\pi}{4}}}{z-e^{-\frac{i\pi}{4}}}
+\frac{3}{4}\frac{1}{z+e^{-\frac{i\pi}{4}}}\right)+{\rm
Re}\left(\frac{\sqrt{2}}{4}\frac{1}{(z+e^{-\frac{i\pi}{4}})^2}+\frac{1}{4}\frac{1}{z+e^{-\frac{i\pi}{4}}}\right) -\frac{\sqrt{2}\pi}{8}-\frac{\sqrt{2}}{2},
$$
where $u={\rm Re\,}(h(z)+g(z))$ and $v={\rm Im\,}(h(z)-g(z))$.

In this case, according to Theorem \Ref{Thm A}, we have
\beqq
\phi_3(z) &=&  2ibh'(z)= \pm 2i\frac{z}{(z-e^{-\frac{i\pi}{4}})(z+e^{-\frac{i\pi}{4}})^3}\\
&=& \pm 2i\left[\frac{i}{8}\frac{1}{z-e^{-\frac{i\pi}{4}}}-\frac{i}{8}\frac{1}{z+e^{-\frac{i\pi}{4}}}
-\frac{1}{4}e^{\frac{i\pi}{4}}\frac{1}{(z+e^{-\frac{i\pi}{4}})^2}+\frac{1}{2}\frac{1}{(z+e^{-\frac{i\pi}{4}})^3}\right].
\eeqq
Integration from $0$ to $z$ gives
$$F(z)=\mp{\rm
Re}\left[\frac{1}{4}\log\frac{z-e^{\frac{-i\pi}{4}}}{z+e^{-\frac{i\pi}{4}}}+\frac{e^{-\frac{i\pi}{4}}}{2}\frac{1}{z+e^{-\frac{i\pi}{4}}}
+\frac{i}{2}\frac{1}{(z+e^{-\frac{i\pi}{4}})^2}\right]+c.
$$

\medskip

\noindent {\bf Case 4:}\quad Let  $\gamma=\frac{7\pi}{4}$.

In this case, $h'(z)$ given by \eqref{eq4b} takes the form
$$h'(z)=\frac{1}{(z+e^{\frac{i\pi}{4}})(z-e^{\frac{i\pi}{4}})^3}
$$
and therefore,
$$h'(z)=
-\frac{1}{8}e^{\frac{i\pi}{4}}\left (\frac{1}{z-e^{\frac{i\pi}{4}}} -\frac{1}{z+e^{\frac{i\pi}{4}}}\right )
+\frac{1}{4}e^{\frac{i\pi}{2}}\frac{1}{(z-e^{\frac{i\pi}{4}})^2}+\frac{1}{2}e^{-\frac{i\pi}{4}}\frac{1}{(z-e^{\frac{i\pi}{4}})^3}.
$$
Integration from $0$ to $z$ gives
\be\label{eq4i}
h(z)=-\frac{1}{8}e^{\frac{i\pi}{4}}\log\frac{z-e^{\frac{i\pi}{4}}}{z+e^{\frac{i\pi}{4}}}-\frac{i}{4}\frac{1}{z-e^{\frac{i\pi}{4}}}
-\frac{e^{-\frac{i\pi}{4}}}{4}\frac{1}{(z-e^{\frac{i\pi}{4}})^2}-\frac{\pi}{8}e^{-\frac{i\pi}{4}}-\frac{1}{2}e^{\frac{i\pi}{4}}.
\ee
Using \eqref{eq4a} for $\gamma=\frac{7\pi}{4}$, we find that
$$h(z)+g(z)=-\frac{1}{z-e^{\frac{i\pi}{4}}} +\sqrt{2}e^{\frac{i\pi}{4}}h(z)- e^{-\frac{i\pi}{4}}
$$
and
$$h(z)-g(z)=\frac{1}{z-e^{\frac{i\pi}{4}}} +\sqrt{2}e^{-\frac{i\pi}{4}}h(z)+ e^{-\frac{i\pi}{4}},
$$
where $h$ is defined by \eqref{eq4i}. We thus obtain that
$$u=-\frac{\sqrt{2}\pi}{8}-\frac{\sqrt{2}}{2}+{\rm
Im}\left(\frac{\sqrt{2}}{8}\log\frac{z-e^{\frac{i\pi}{4}}}{z+e^{\frac{i\pi}{4}}}
+\frac{1}{4}\frac{1}{z-e^{\frac{i\pi}{4}}}\right)-{\rm
Re}\left(\frac{\sqrt{2}}{4}\frac{1}{(z-e^{\frac{i\pi}{4}})^2}+\frac{3}{4}\frac{1}{z-e^{\frac{i\pi}{4}}}\right),
$$
and
$$v=\frac{\sqrt{2}\pi}{8}-\frac{\sqrt{2}}{2}-{\rm
Im}\left(\frac{\sqrt{2}}{8}\log\frac{z-e^{\frac{i\pi}{4}}}{z+e^{\frac{i\pi}{4}}}
-\frac{3}{4}\frac{1}{z-e^{\frac{i\pi}{4}}}\right)+{\rm
Re}\left(\frac{\sqrt{2}}{4}\frac{1}{(z-e^{\frac{i\pi}{4}})^2}-\frac{1}{4}\frac{1}{z-e^{\frac{i\pi}{4}}}\right),
$$
where $u={\rm Re\,}(h(z)+g(z))$ and $v={\rm Im\,}(h(z)-g(z))$.

In this case, by Theorem \Ref{Thm A}, we find that
\beqq
\phi_3(z) &=& 2ibh'(z)=\pm 2i\frac{z}{(z+e^{\frac{i\pi}{4}})(z-e^{\frac{i\pi}{4}})^3}\\
&=&\pm 2i\left[-\frac{i}{8}\frac{1}{z+e^{\frac{i\pi}{4}}}+\frac{i}{8}\frac{1}{z-e^{\frac{i\pi}{4}}}
+\frac{1}{4}e^{-\frac{i\pi}{4}}\frac{1}{(z-e^{\frac{i\pi}{4}})^2}+\frac{1}{2}\frac{1}{(z-e^{\frac{i\pi}{4}})^3}\right].
\eeqq
Integration from $0$ to $z$ gives
$$F(z)=\mp {\rm
Re}
\left[-\frac{1}{4}\log\frac{z+e^{\frac{i\pi}{4}}}{z-e^{\frac{i\pi}{4}}}+\frac{1}{2}e^{\frac{i\pi}{4}}\frac{1}{z-e^{\frac{i\pi}{4}}}
+\frac{i}{2}\frac{1}{(z-e^{\frac{i\pi}{4}})^2}\right]+c.
$$
\medskip

\noindent {\bf Case 5:}\quad Let  $\gamma\notin\Big\{\frac{\pi}{4},\
\frac{3\pi}{4},\ \frac{5\pi}{4},\ \frac{7\pi}{4}\Big\}.$

In this case, $h'(z)$ given by \eqref{eq4b}  has simple poles at $ie^{i\gamma}$ and $-ie^{i\gamma}$, and a pole of
order $2$ at $e^{-i\gamma}$. Thus, we may rewrite $h'(z)$ as
$$h'(z)= \frac{A}{z+ie^{i\gamma}} + \frac{B}{z-ie^{i\gamma}} +\frac{C}{z-e^{-i\gamma}}
+\frac{D}{(z-e^{-i\gamma})^2}
$$
where $A,B,C$ and $ D$ can be easily computed using a standard procedure from residue calculus or otherwise. Indeed
$$A=\frac{e^{-i\gamma}}{4(1-\sin 2\gamma)}, ~B=\frac{e^{-i\gamma}}{4(1+\sin 2\gamma)}, ~C=-\frac{e^{-i\gamma}}{2\cos^2 2\gamma},~\mbox{ and }~
D=\frac{1}{2\cos 2\gamma}.
$$
We observe that $A+B+C=0$.
Integration from $0$ to $z$ leads to
\be\label{eq4j}
h(z)= A\log\frac{z+ie^{i\gamma}}{ie^{i\gamma}}+ B\log\frac{z-ie^{i\gamma}}{-ie^{i\gamma}}
+C\log\frac{z-e^{-i\gamma}}{-e^{-i\gamma}}-
 \frac{D}{z-e^{-i\gamma}}- De^{i\gamma} .
\ee
Note that $g$ defined by \eqref{eq4a} gives
$$
h(z)+g(z)= -\frac{1}{z-e^{-i\gamma}} -2i e^{i\gamma}\sin \gamma h(z)-e^{i\gamma}
$$
and
$$h(z)-g(z)=\frac{1}{z-e^{-i\gamma}} +2e^{i\gamma}\cos \gamma  h(z)+ e^{i\gamma}
$$
where $h$ is defined by \eqref{eq4j}.
By computation, we know that
$$u={\rm
Im}\left(\frac{\sin\gamma}{2(1-\sin
2\gamma)}\log\frac{z+ie^{i\gamma}}{ie^{i\gamma}}+\frac{\sin\gamma}{2(1+\sin
2\gamma)}\log\frac{z-ie^{i\gamma}}{-ie^{i\gamma}}-\frac{\sin\gamma}{\cos^2
2\gamma}\log\frac{z-e^{-i\gamma}}{-e^{-i\gamma}}\right)
$$
$$-\frac{\cos\gamma}{\cos2\gamma}{\rm Re}\left (\frac{z}{z-e^{-i\gamma}}\right ),
$$
$$v={\rm Im}\left(\frac{\cos\gamma}{2(1-\sin
2\gamma)}\log\frac{z+ie^{i\gamma}}{ie^{i\gamma}}+\frac{\cos\gamma}{2(1+\sin
2\gamma)}\log\frac{z-ie^{i\gamma}}{-ie^{i\gamma}}-\frac{\cos\gamma}{\cos^2
2\gamma}\log\frac{z-e^{-i\gamma}}{-e^{-i\gamma}}\right)
$$
$$-\frac{\sin\gamma}{\cos2\gamma}{\rm Re}\left (\frac{z}{z-e^{-i\gamma}}\right ) .
$$


In the final case, by Theorem \Ref{Thm A}, we find that
\beqq
\phi_3(z) &=&  2ibh'(z)= \pm \frac{2iz}{(z^2+e^{2i\gamma})(z-e^{-i\gamma})^2}\\
&=& \pm2i\left[-\frac{i}{4(1-\sin 2\gamma)(z+ie^{i\gamma})}+\frac{i}{4(1+\sin 2\gamma)(z-ie^{i\gamma})} \right. \\
&&  \hspace{1cm}+ \left .\frac{i\sin 2\gamma}{2(z-e^{-i\gamma})\cos^2
2\gamma}+\frac{e^{-i\gamma}}{2(z-e^{-i\gamma})^2\cos
2\gamma}\right].
\eeqq
Integration from $0$ to $z$ gives
$$F=\pm{\rm
Re}\left[\frac{\log(z+ie^{i\gamma})}{2(1-\sin
2\gamma)}-\frac{\log(z-ie^{i\gamma})}{2(1+\sin 2\gamma)}-\frac{\sin
2\gamma}{\cos^2
2\gamma}\log(z-e^{-i\gamma})-\frac{ie^{-i\gamma}}{(z-e^{-i\gamma})\cos
2\gamma}\right]+c.
$$
The proof is complete. \epf

\section{Vertical strips}\label{strips}

Hengartner and Schober \cite{HS87} investigated the family of functions from ${\mathcal S}_H$ that map $\ID$
onto the horizontal strip domain $\{w:\,| {\rm Im}\, w| < \pi/4\}$. As an analogous result,  Dorff \cite{Do99}
considered the family ${\mathcal S}_H(\ID, \Omega_\alpha)$ of functions from  ${\mathcal S}_H$ which map
$\ID$ onto the asymmetric vertical strip domains
$$\Omega_\alpha=\Big \{w:\, \frac{\alpha-\pi}{2\sin\alpha}<{\rm Re\,}w<\frac{\alpha}{2\sin\alpha}\Big \},
$$
where $\frac{\pi}{2}\leq\alpha<\pi$. Set ${\mathcal S}_H^0(\ID, \Omega_\alpha)={\mathcal S}_H(\ID, \Omega_\alpha)\cap {\mathcal S}_H^0.$
Note that $\Omega_{\pi /2}=\{w: \,| {\rm Re}\, w| < \pi/4\}$ and so, the class discussed by Hengartner and Schober \cite{HS87}
follows by using a suitable rotation.

\begin{lem}\label{dila2}
Each $f=h+\overline{g} \in {\mathcal S}_H^0(\ID, \Omega_\alpha)$  has the form
\be\label{eq5}
h(z)+g(z)=\psi (z), \quad \psi(z)= \frac{1}{2i\sin\alpha}\log\left(\frac{1+ze^{i\alpha}}{1+ze^{-i\alpha}}
\right).
\ee
Moreover,
\be\label{eq5a}
h'(z)=\frac{\psi'(z)}{1+\omega(z)},~
g'(z)=\frac{\omega(z) \psi'(z)}{1+\omega(z)} ~\mbox{ and }~ \psi'(z)= \frac{1}{(1+ze^{-i\alpha})(1+ze^{i\alpha})}.
\ee
Here $\omega (z)=g'(z)/h'(z)$ denotes the dilatation of $f$.
\end{lem}
\bpf
The representation \eqref{eq5} is well-known whereas \eqref{eq5a} follows if we solve the pair of equations
$h'(z)+g'(z)=\psi '(z)$ and $\omega (z)h'(z)-g'(z)=0$.
The proof is complete. \epf

\begin{thm}\label{th2}
Let $S$ be a minimal surface over $\Omega_\alpha$ with the
projection $f=h+\overline{g}\in {\mathcal S}_H^0(\ID, \Omega_\alpha)$,
which satisfies \eqref{eq2} and whose dilatation $\omega=b^2$, where
$b(z)=\pm z$. Then $S=\{(u, v, F(u,v) ):\, u+iv\in\Omega_\alpha\},$ where
\beqq
u& =& \ds \frac{1}{2\sin\alpha}{\rm Im}\left [\log \left (\frac{1+ze^{i\alpha}}{1+ze^{-i\alpha}}\right)\right ],\\
v &= & \left\{
\begin{array}{ll}
\ds {\rm Im}\left (\frac{z}{z^2+1} \right ) & \mbox{ if }~ \alpha=\frac{\pi}{2}, \\
\ds \frac{1}{2\cos\alpha}{\rm Im}\left[\log \left (\frac{(1+ze^{i\alpha})(1+ze^{-i\alpha})}{z^2+1}\right) \right]
&  \mbox{ if }~ \frac{\pi}{2}<\alpha<\pi
\end{array}
\right.
\eeqq
and
$$F=\left\{
\begin{array}{ll}
\ds \pm {\rm Im} \left (\frac{1}{z^2+1}\right )+c & \mbox{ if }~\alpha=\frac{\pi}{2}, \\
\ds \pm{\rm Re}\left [\frac{1}{2\cos\alpha}\log\left(\frac{z+i}{z-i}\right )-\frac{1}{\sin{2\alpha}}\log\left (\frac{z+e^{i\alpha}}{z+e^{-i\alpha}}\right)\right]+c
& \mbox{ if }~\frac{\pi}{2}<\alpha<\pi.
\end{array}
\right.
$$
\end{thm}
\bpf  Let $f=h+\overline{g}\in {\mathcal S}_H^0(\ID, \Omega_\alpha)$ with $\omega (z)=z^2$. Then by Lemma \ref{dila2}, we have
\be\label{eq5d}
h'(z) = \left\{
\begin{array}{ll}
\ds  \frac{1}{(z+i)^2(z-i)^2} & \mbox{ if }~\alpha=\frac{\pi}{2},\\
\ds  \frac{1}{(z+i)(z-i)(z+e^{i\alpha})(z+e^{-i\alpha})}& \mbox{ if }~\frac{\pi}{2}<\alpha<\pi .
\end{array}
\right.
\ee
{\bf Case (i):} Let $\alpha=\frac{\pi}{2}$. Then consider the partial fraction expression for $h'(z)$:
$$h'(z)=\frac{1}{4}\left[\frac{i}{z+i}-\frac{i}{z-i}-\frac{1}{(z+i)^2}-\frac{1}{(z-i)^2}\right].
$$
Integration from $0$ to $z$ gives
\be\label{eq5b}
h(z)=\frac{1}{4}\left[i\log (z+i)-i\log (z-i)+\frac{1}{z+i}+\frac{1}{z-i}+\pi\right].
\ee
Also, by \eqref{eq5}, we obtain that
\be\label{eq5c}
h(z)+g(z)= \frac{1}{2i}\log\left(\frac{i-z}{i+z}\right)=\frac{1}{2}\left[i\log (z+i)-i\log (z-i)+\pi\right]
\ee
so that, by \eqref{eq5b} and \eqref{eq5c}
$$h(z)-g(z)= 2h(z)-(h(z)+g(z))=  \frac{z}{z^2+1} .
$$
As before, it follows that
$$u=\frac{1}{2}{\rm Im}\left(\log\left(\frac{i-z}{i+z}\right)\right ), ~
v={\rm Im} \left(\frac{z}{z^2+1}\right )
$$
and $\phi_3$ given by Theorem \Ref{Thm A} takes the form
$$\phi_{3}(z)=\pm\frac{2iz}{(1-iz)^2(1+iz)^2}=\pm \frac{1}{2}\left (\frac{1}{(z-i)^2}-\frac{1}{(z+i)^2} \right ).
$$
We thus obtain $F$  by integration:
$$F=\pm {\rm Im}\left (\frac{1}{z^2+1}\right ) +c.
$$

\medskip
\noindent
{\bf Case (ii):} Let $\frac{\pi}{2}<\alpha<\pi$. The partial fraction expansion of $h'(z)$ in \eqref{eq5d} yields
$$h'(z)= -\frac{1}{4\cos\alpha}\left(\frac{1}{z+i}+\frac{1}{z-i}\right)+\frac{1}{(e^{-i\alpha}-e^{3i\alpha})(z+e^{i\alpha})}
+\frac{1}{(e^{i\alpha}-e^{-3i\alpha})(z+e^{-i\alpha})}.
$$
Integration from $0$ to $z$ gives
$$h(z)=-\frac{1}{4\cos\alpha}\log(z^2+1)+\frac{1}{e^{-i\alpha}-e^{3i\alpha}}\log(1+ze^{-i\alpha})
+\frac{1}{e^{i\alpha}-e^{-3i\alpha}}\log(1+ze^{i\alpha}),
$$
which simplifies to
\be\label{eq5e}
h(z)=-\frac{1}{4\cos\alpha}\log(z^2+1)+\frac{ie^{-i\alpha}}{2\sin 2\alpha}\log(1+ze^{-i\alpha})
-\frac{ie^{i\alpha}}{2\sin 2\alpha}\log(1+ze^{i\alpha}).
\ee
By using \eqref{eq5}, we obtain that
$$u= {\rm Re}\, (h(z)+g(z)) = \frac{1}{2\sin\alpha}{\rm Im}\left(\log\left(\frac{1+ze^{i\alpha}}{1+ze^{-i\alpha}}\right )\right).
$$
Writing $h(z)-g(z)= 2h(z)-(h(z)+g(z))$ and using \eqref{eq5} and \eqref{eq5e}, we can easily find that
$$h(z)-g(z)=-\frac{1}{2\cos\alpha}\log(z^2+1) +\frac{1}{2\cos\alpha}\log(1+ze^{-i\alpha})
+\frac{1}{2\cos\alpha}\log(1+ze^{i\alpha})
$$
which gives
$$v={\rm Im}\, (h(z)-g(z)) =\frac{1}{2\cos\alpha}{\rm Im}\left(\log\frac{(1+ze^{i\alpha})(1+ze^{-i\alpha})}{z^2+1}\right).
$$
In this case, $\phi_3$ given by Theorem \Ref{Thm A} takes the form
\beqq
\phi_{3}(z)& =&\pm\frac{2iz}{(z+i)(z-i)(z+e^{i\alpha})(z+e^{-i\alpha})}\\
&= &\pm 2i\left[\frac{i}{4\cos\alpha}\left(\frac{1}{z+i}-\frac{1}{z-i}\right)
+\frac{1}{2i\sin{2\alpha}}\left(\frac{1}{z+e^{i\alpha}}-\frac{1}{z+e^{-i\alpha}}\right)\right].
\eeqq
Integration from $0$ to $z$ gives
$$F=\pm{\rm
Re}\left[\frac{1}{2\cos\alpha}\log\left(\frac{z+i}{z-i}\right)-\frac{1}{\sin{2\alpha}}\log\left(\frac{z+e^{i\alpha}}{z+e^{-i\alpha}}\right)\right]+c
$$
and the proof is complete. \epf

\section{Single slit}

Finally, we consider single slit domain $L$ whose slit lies on the negative real
axis.
Moreover, by the result of Livingston \cite{Living92} (see also \cite{Living97} and  Dorff \cite[Corollary 2]{Do99})
it follows that if $f=h+\overline{g} \in
{\mathcal S}_H^0$ is a slit mapping whose slit lies on the negative
real axis, then one has
\be\label{eq6}
h(z)-g(z)=\frac{z}{(1-z)^2}.
\ee

\begin{thm}\label{th3}
Let $S$ be a minimal surface over $L$ with the projection $f=h+\overline{g}\in {\mathcal S}_H^0$, which satisfies
\eqref{eq6} and whose dilatation $\omega=b^2$, where $b(z)=\pm z$. Then $S=\{(u, v, F(u,v) ):\, u+iv\in L\},$
where
$$u={\rm Re}\left(\frac{2z^3-3z^2+3z}{3(1-z)^3}\right), ~
v={\rm Im}\left(\frac{z}{(1-z)^2}\right),~
$$
and
$$F=\pm{\rm Im}\left(\frac{1}{(z-1)^2}+\frac{2}{3(z-1)^3}\right)+c.
$$
\end{thm}\bpf
By assumption, $f=h+\overline{g}\in {\mathcal S}_H^0$ is a single slit mapping whose slit lies on the
negative real axis with $\omega (z)=z^2$. Then \eqref{eq6} holds and therefore, we have
$$h'(z)-g'(z)=\frac{1+z}{(1-z)^3} ~\mbox{ and }~\ g'(z)=z^2h'(z).
$$
Solving these two equations, we obtain
$$h'(z)=\frac{1}{(1-z)^4} 
.
$$
Integrating from $0$ to $z$ yields
$$h(z)=-\frac{1}{3}+\frac{1}{3(1-z)^3}
$$
and so
$$ g(z)=h(z)- \frac{z}{(1-z)^2}= -\frac{1}{3}+\frac{1}{3(1-z)^3} -\frac{z}{(1-z)^2},
$$
which, by using the previous equation, gives
$$h(z)+g(z)=\frac{2z^3-3z^2+3z}{3(1-z)^3}.
$$
The desired representations for $u= {\rm Re} \,(h(z)+g(z))$ and $v= {\rm Im} \,(h(z)-g(z))$ follow easily.
Finally, since
$$ \phi_{3}(z)=\pm 2izh'(z)=\pm\frac{2iz}{(1-z)^4}=\pm 2i\left (\frac{1}{(1-z)^4}-\frac{1}{(1-z)^3}\right ),
$$
integrating this from $0$ to $z$ yields
$$F=\pm{\rm Im}\left(\frac{1}{(z-1)^2}+\frac{2}{3(z-1)^3}\right)+c.
$$
The proof is complete.
 \epf

\section{Illustration using Mathematica}

\begin{figure}[htp]
\begin{center}
\includegraphics[height=6cm, width=5.5cm, scale=4.5]{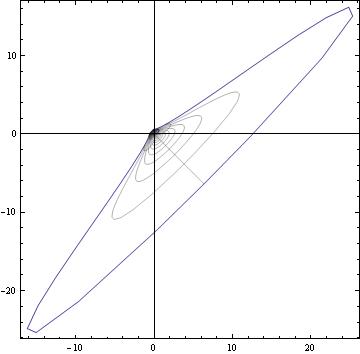}
\hspace{0.5cm}
\includegraphics[height=6cm, width=5.5cm, scale=4.5]{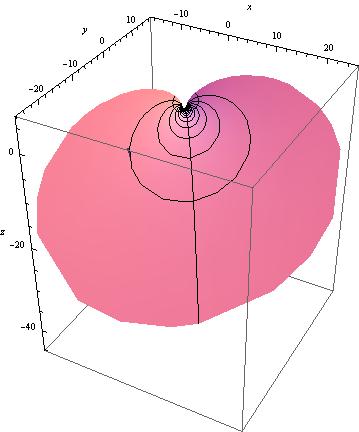}
\end{center}
\caption{\textbf{Case 1:} $\gamma =\pi /4$ of Theorem \ref{1} \label{Th1Case1:fig1}}
\end{figure}

\begin{figure}[htp]
\begin{center}
\includegraphics[height=6cm, width=5.5cm, scale=4.5]{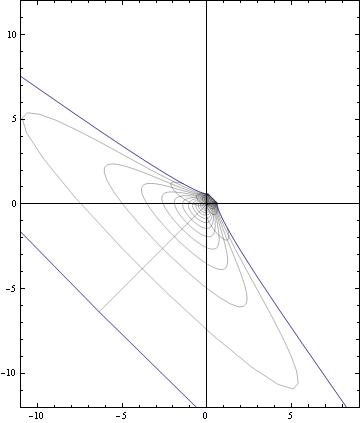}
\hspace{0.5cm}
\includegraphics[height=6cm, width=5.5cm, scale=4.5]{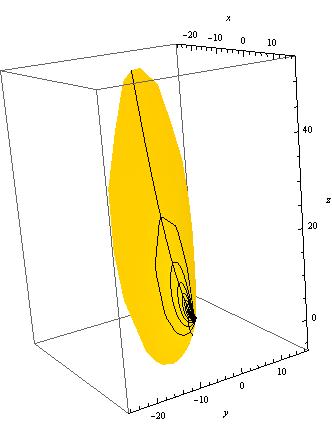}
\end{center}
\caption{\textbf{Case 2:} $\gamma =3\pi /4$ of Theorem \ref{1} \label{Th1Case2:fig2}}
\end{figure}

\begin{figure}[htp]
\begin{center}
\includegraphics[height=6cm, width=5.5cm, scale=4.5]{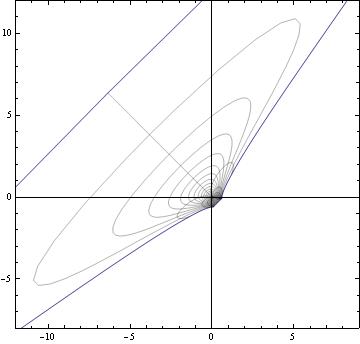}
\hspace{0.6cm}
\includegraphics[height=6cm, width=5.5cm, scale=20]{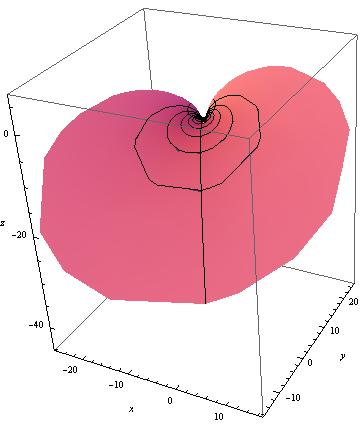}
\end{center}
\caption{\textbf{Case 3:} $\gamma =5\pi /4$ of Theorem \ref{1} \label{Th1Case3:fig3}}
\end{figure}

\begin{figure}[htp]
\begin{center}
\includegraphics[height=6cm, width=5.5cm, scale=4.5]{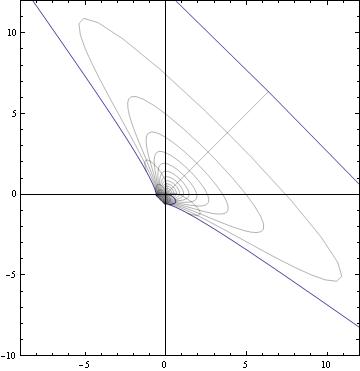}
\hspace{0.6cm}
\includegraphics[height=6cm, width=5.5cm, scale=9.5]{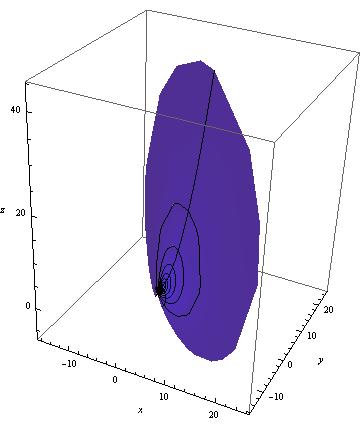}
\end{center}
\caption{\textbf{Case 4:} $\gamma =7\pi /4$ of Theorem \ref{1} \label{Th1Case5:fig4}}
\end{figure}

\begin{figure}[htp]
\begin{center}
\includegraphics[height=6cm, width=5.5cm, scale=4.5]{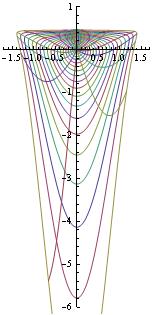}
\hspace{0.6cm}
\includegraphics[height=6cm, width=5.5cm, scale=9.5]{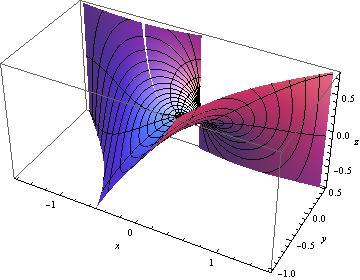}
\end{center}
\caption{\textbf{Case 5} with $\gamma =\pi /2$ of Theorem \ref{1} \label{Th1Case5:fig5}}
\end{figure}

\begin{figure}[htp]
\begin{center}
\includegraphics[height=6cm, width=5.5cm, scale=4.5]{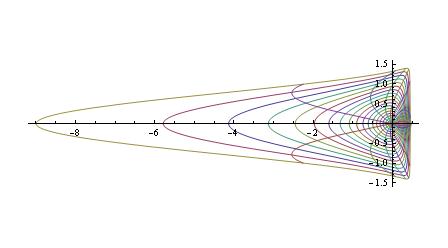}
\hspace{0.6cm}
\includegraphics[height=6cm, width=5.5cm, scale=9.5]{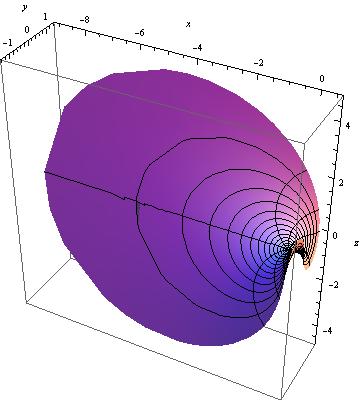}
\end{center}
\caption{\textbf{Case 5} with $\gamma =\pi$ of Theorem \ref{1} \label{Th1Case5:fig6}}
\end{figure}

The images of the disk $|z|<r$ for $r$ closer to $1$ under $f=h+\overline{g}$ for various cases of
Theorem \ref{1} and the corresponding
minimal surfaces associated with $f$ are illustrated in Figures \ref{Th1Case1:fig1}-\ref{Th1Case5:fig7}.
Similar illustrations for Theorem \ref{th2} (see Figures \ref{Th2fig1}-\ref{Th2fig4}) and
Theorem \ref{th3} (see Figure \ref{Th3fig1}) are also provided. These figures are drawn
using Mathematica (see for example \cite{Ruskeepaa}).

\begin{figure}[htp]
\begin{center}
\includegraphics[height=6cm, width=5.5cm, scale=4.5]{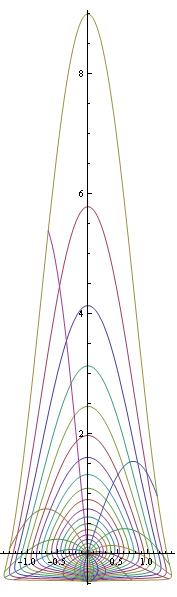}
\hspace{0.6cm}
\includegraphics[height=6cm, width=5.5cm, scale=9.5]{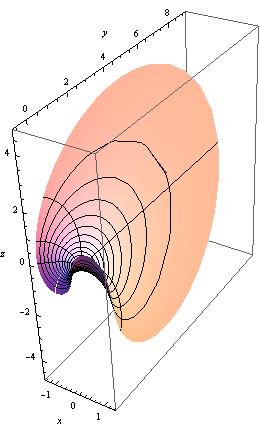}
\end{center}
\caption{\textbf{Case 5} with $\gamma =3\pi/2$ of Theorem \ref{1}}
\end{figure}

\begin{figure}[htp]
\begin{center}
\includegraphics[height=6cm, width=5.5cm, scale=4.5]{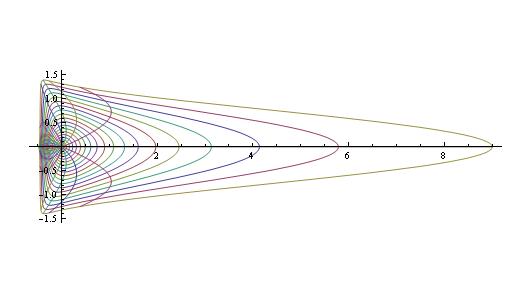}
\hspace{0.6cm}
\includegraphics[height=6cm, width=5.5cm, scale=9.5]{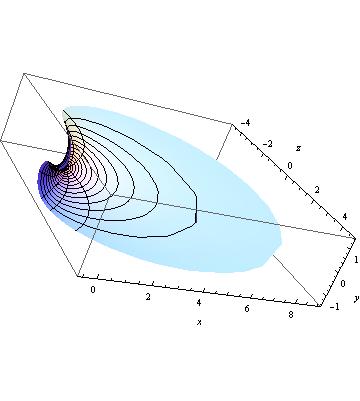}
\end{center}
\caption{\textbf{Case 5} with $\gamma =0$ of Theorem \ref{1} \label{Th1Case5:fig7}}
\end{figure}

\begin{figure}[htp]
\begin{center}
\includegraphics[height=6cm, width=5.5cm, scale=4.5]{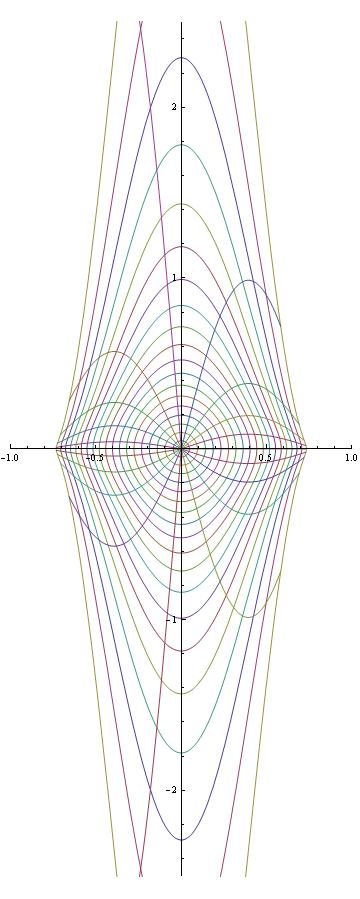}
\hspace{0.6cm}
\includegraphics[height=6cm, width=5.5cm, scale=9.5]{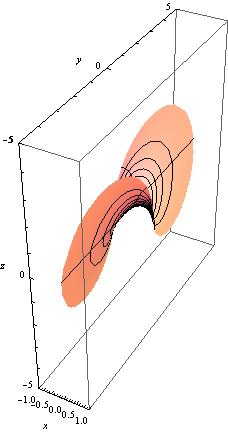}
\end{center}
\caption{Illustration for $\alpha =\pi/2$ of Theorem \ref{th2} \label{Th2fig1}}
\end{figure}


\begin{figure}[htp]
\begin{center}
\includegraphics[height=6cm, width=5.5cm, scale=4.5]{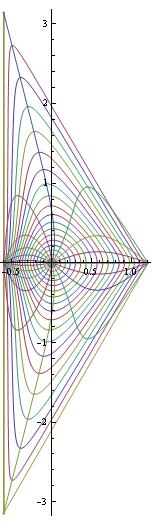}
\hspace{0.6cm}
\includegraphics[height=6cm, width=5.5cm, scale=9.5]{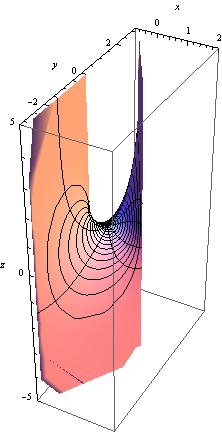}
\end{center}
\caption{Illustration for $\alpha =2\pi/3$ of Theorem \ref{th2} \label{Th2fig2}}
\end{figure}

\begin{figure}[htp]
\begin{center}
\includegraphics[height=6cm, width=5.5cm, scale=4.5]{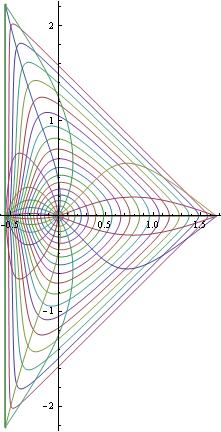}
\hspace{0.6cm}
\includegraphics[height=6cm, width=5.5cm, scale=9.5]{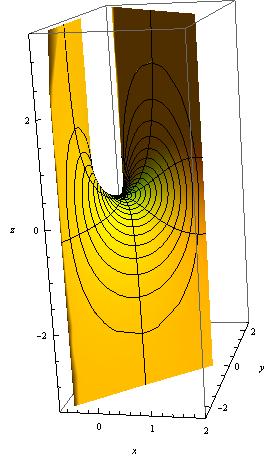}
\end{center}
\caption{Illustration for $\alpha =3\pi/4$ of Theorem \ref{th2} \label{Th2fig3}}
\end{figure}

\begin{figure}[htp]
\begin{center}
\includegraphics[height=6cm, width=5.5cm, scale=4.5]{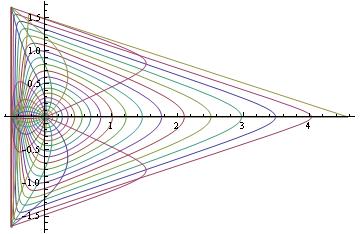}
\hspace{0.6cm}
\includegraphics[height=6cm, width=5.5cm, scale=9.5]{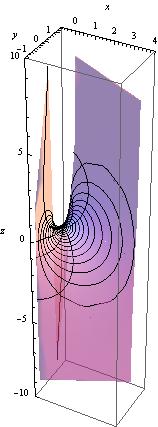}
\end{center}
\caption{Illustration for $\alpha =9\pi/10$ of Theorem \ref{th2} \label{Th2fig4}}
\end{figure}

\begin{figure}[htp]
\begin{center}
\includegraphics[height=6cm, width=5.5cm, scale=4.5]{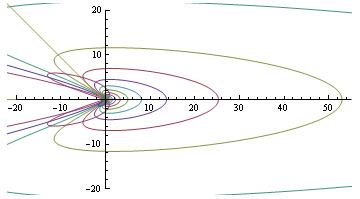}
\hspace{0.6cm}
\includegraphics[height=6cm, width=5.5cm, scale=9.5]{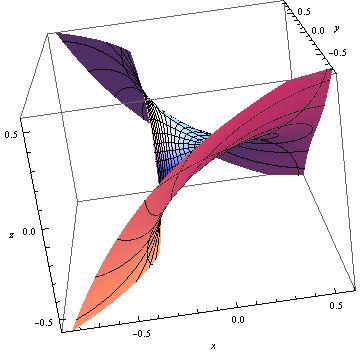}
\end{center}
\caption{Illustration for Theorem \ref{th3} \label{Th3fig1}}
\end{figure}


\end{document}